\magnification 1200 \baselineskip=16pt
\input amssym

\def\n{\noindent}
\def\qed{\qquad {\rm qed}}
\def\qed{{\hfill{\vrule height7pt width7pt depth0pt}\par\bigskip}}

\def\8{\infty}

\def\M{{\cal M}}

\def\E{{\cal E}}

\centerline{\bf On the best constants in some non-commutative
martingale inequalities}
\bigskip

\setbox1=\vtop{\hsize=6cm \centerline {Marius Junge$^{\dag}$}
\centerline{University of Illinois} \centerline{Department of
Mathematics} \centerline{Urbana, IL 61801} \centerline{USA}
\centerline{junge@math.uiuc.edu}}

\setbox2=\vtop{\hsize=6cm \centerline{Quanhua Xu}
\centerline{Universit\'e de Franche-Comt\'e}
\centerline{Laboratoire de Math\'ematiques} \centerline{25030
Besan\c con Cedex} \centerline{France}
\centerline{qx@math.univ-fcomte.fr}}
$$\line{\box1\hfill\box2}$$

\vskip 1.5cm

\noindent{\bf Abstract.} We determine the optimal orders for the
best constants in the non-commutative Burkholder-Gundy, Doob  and
Stein inequalities obtained recently in the non-commutative
martingale theory.

\vfill



 \footnote{}{$^{\dag}$Marius
Junge is partially supported by the NSF}

\vfill\eject

 \centerline{\bf 1. Introduction}
\medskip

\n The theory of non-commutative martingale inequalities has been
rapidly developed for several years. Many of the classical
inequalities in the usual martingale theory have been already
transferred into the non-commutative setting. As in  the
commutative case (see [B1, B2] and the references therein) the
order of the best constants in these inequalities may provide
important additional insight. We refer to [Ru] for applications of
the order in the non-commutative Khintchine inequalities and to
the work of Nazarov, Pisier, Treil and  Volberg  [NPTV] for
applications of the order of the UMD-constant for Schatten
classes. In this note we analyze the order of the best constants
for the non-commutative Burkholder-Gundy inequalities and the
non-commutative Stein inequalities proved in [PX1-2] as well as
the non-commutative Doob maximal inequalities in [J]. We refer to
[JX1] for the known constants in the Burkholder and Rosenthal
inequalities.

Note that the non-commutative Burkholder-Gundy inequalities imply
that the non-commutative martingale transforms by sequences of
signs are all of type $(p,p)$ for every $1<p<\8$. However, the
problem whether these non-commutative martingale transforms are of
weak type $(1,1)$ was left open since [PX1-2]. Only very recently
that this problem was affirmatively solved by Randrianantoanina
[R]. This yields the optimal order of the unconditionality
constant for martingale differences. This improves considerably on
the exponential estimates in [PX1-2], [J] and [JX1]. The purpose
of this note is to further clarify the optimal order of the best
constants in the martingale inequalities mentioned above. We
obtain the optimal order except for one case (see below).

Let us recall the notion of non-commutative martingales and the
formulation of the martingale inequalities considered in the note.
Throughout this text ${\cal M}$  denotes a finite von Neumann
algebra equipped with a normal faithful normalized trace $\tau$,
and $({\cal M}_{n})_{n \geq 1}$ an increasing filtration of von
Neumann subalgebras of ${\cal M}$ whose union is w*-dense in
${\cal M}.$ For $1 \leq p \leq \infty$ we denote by $L^{p} ({\cal
M}, \tau),$ or simply $L^{p}({\cal M})$ the usual non-commutative
$L^{p}$-space associated with $({\cal M}, \tau)$ (cf., e.g. [D],
[S]; see also the survey [PX3]).  Recall that by convention,
$L^{\infty}({\cal M})={\cal M}$ with the operator norm. As usual,
$L^{p}({\cal M}_{n}) = L^{p}({\cal M}_{n}, \tau \big\vert_{{\cal
M}_{n}})$ is naturally identified as a subspace of $L^{p}({\cal
M}).$ It is well-known that there is a unique normal faithful
conditional expectation ${\cal E}_{n}$ from ${\cal M}$ onto ${\cal
M}_{n}$ such that $\tau \circ {\cal E}_{n} = \tau.$ Moreover,
${\cal E}_{n}$ extends to a contractive projection from
$L^{p}({\cal M})$ onto $L^{p}({\cal M}_{n}),$ for every $1\le
p<\8$, which is still denoted by ${\cal E}_{n}.$

A non-commutative martingale with respect to $({\cal M}_{n})$ is a
sequence $x = (x_{n})_{n \geq 1}$ in $ L^{1}({\cal M})$ such that
$$x_{n} = {\cal E}_{n} (x_{n+1}),\quad \forall\ n \geq 1.$$
The difference sequence of $x$ is $d x = (d x_{n})_{n \geq 1},$
where $d x_{n} = x_{n} - x_{n-1}$ (with $x_{0} = 0$ by
convention). Then we define $L^{p}$-martingales and bounded
$L^{p}$-martingales, as usual. If $x$ is an $L^{p}$-martingale, we
set
$$\Vert x \Vert_{p} = \sup_{n} \, \Vert x_{n} \Vert_{p}\,.$$
In the sequel, we will fix ${\cal M}$, $\tau$ and $({\cal
M}_{n})_{n \geq 1}$ as above. all martingales will be
non-commutative martingales with respect to the fixed filtration
$({\cal M}_{n})_{n \geq 1}$, unless explicitly stated otherwise.

\smallskip

To state the non-commutative Burkholder-Gundy inequalities, we
introduce the norms in the Hardy spaces of martingales defined in
[PX1-2]. Let $1 \leq p \leq \infty$ and $x = (x_{n})_{n \geq 1}$
be an $L^{p}$-martingale. Set, for $p \geq 2$
$$\Vert x \Vert_{{\cal H}^{^p}} = \ \max\ \biggl\{ \bigl\Vert \bigl(
\sum_{n \geq 1} \vert d x_{n} \vert^{2}\bigl)^{1/2} \bigl\Vert_{p},\
\bigl\Vert \bigl( \sum_{n \geq 1} \vert d x^{\ast}_{n}\vert^{2}
\bigl)^{1/2} \bigl\Vert_{p} \biggl\}$$
and for $p < 2$
$$\Vert x \Vert_{{\cal H}^{p}} =\  \inf\ \biggl\{ \bigl\Vert \bigl(
\sum_{n \geq 1} \vert d y_{n}\vert^{2} \bigl)^{1/2} \bigl\Vert_{p}
+ \bigl\Vert \bigl( \sum_{n \geq 1} \vert d z^{\ast}_{n}\vert^{2}
\bigl)^{1/2} \bigl\Vert_{p} \biggl\},$$ where the infimum runs
over all decompositions $x = y + z$ of $x$ as sums of two
$L^{p}$-martingales. Recall that $|\,\cdot\,|$ stands for the
usual (right) modulus of operators, i.e. $|a|=(a^*a)^{1/2}$.

\smallskip

 Then the mentioned
non-commutative Burkholder-Gundy inequalities read as follows. In
all what follows, letters $\alpha_{p}, \beta_{p}$, etc $\ldots$
will denote positive constants depending only on $p$, and $C$ an
absolute positive constant.

\proclaim Theorem 1 (Non-commutative Burkholder-Gundy
inequalities). Let $1 < p < \infty.$ Then for all finite
non-commutative $L^{p}$-martingales $x = (x_{n})$
$$\alpha^{-1}_{p} \Vert x \Vert_{{\cal H}^{p}} \leq \Vert x \Vert_{p}
\leq \beta_{p} \ \Vert x \Vert_{{\cal H}^{p}}. \leqno(BG_{p})$$

\smallskip

Inequalities $(BG_{p})$ were first proved in [PX]; see also [R]
for another proof. For Clifford martingales, some particular cases
of $(BG_{p})$ also appear in [CK]. Since the norm $\Vert\, \cdot\,
\Vert_{{\cal H}^{p}}$ is unconditional on martingale difference
sequences, $(BG_{p})$ immediately implies that non-commutative
martingale transforms by sequences of signs are bounded, namely,
we have the following statement.

\proclaim Corollary 2 (Non-commutative martingale transforms: type
$(p,p)$). Let $1 < p < \infty.$ Then for all finite
non-com\-muta\-tive $L^{p}$-martingales
$$\Bigl\Vert \sum_{n \geq 1} \varepsilon_{n} d x_{n}\Bigl\Vert_{p}
\leq \kappa_{p}\ \Vert x \Vert_{p},\quad \forall\ \varepsilon_{n}
= \pm 1. \leqno(MT_{p})$$

\smallskip
 Conversely, by virtue of the non-commutative Khintchine
inequalities (cf. [LLP]), $(M T_{p})$ implies $(BG_{p})$ in the
case of $p > 2$. For the other values of $p,$ we additionally need
the following non-commutative Stein inequality, proved in [PX1]
too. The reader is referred to [R] for another proof of $(S_{p}).$

\smallskip

 \proclaim Theorem 3 (Non-commutative Stein inequality).
Let $1 < p < \infty.$ Then for all finite sequences $(a_{n})_{n
\geq 1}$ in  $L^{p}({\cal M})$
$$\Bigl\Vert \Bigl( \sum\ \vert \E_{n}
a_{n}\vert^{2}\Bigl)^{1/2} \Bigl\Vert_{p} \leq \gamma_{p}\
\Bigl\Vert \Bigl( \sum_{n}\ \vert a_{n}\vert^{2} \Bigl)^{1/2}
\Bigl\Vert_{p}. \leqno(S_{p})$$

\smallskip

The  non-commutative Stein inequality is closely related to
$(BG_p)$ as well as to  the non-commutative Doob inequality
obtained in [J].

\proclaim Theorem 4 (Non-commutative Doob inequality). Let $1 < p
\leq \infty.$ Then for any $a \in L^{p}({\cal M})$ with $a \geq 0$
there is $b \in L^{p} ({\cal M})$ with $b \geq 0$ such that
$$\Vert b \Vert_{p} \leq \delta_{p}\, \Vert a \Vert_{p}\ \hbox{
and }\ \E_{n} a \leq b,\ \forall\ n \geq 0. \leqno(D_{p})$$

Note that in the commutative case the above statement is, of
course, equivalent to the usual Doob inequality on the maximal
functions of martingales. However, in the non-commutative setting
it is unclear how to define the maximal function as an operator.
We refer to [J] for a substitute for the usual maximal function in
the non-commutative case. It is sometimes more convenient to work
with the following dual reformulation of $(D_{p}).$

\proclaim Theorem 4$'$ (Dual form of the non-commutative Doob
inequality). Let $1 \leq p < \infty.$ Then for all finite
sequences $(a_{n})_{n \geq 1}$ of positive elements in
$L^{p}({\cal M})$
$$\Bigl\Vert \sum_{n \geq 1}\ {\cal E}_{n} a_{n} \Bigl\Vert_{p} \leq
\delta'_{p} \Bigl\Vert \sum_{n \geq 1} a_{n}\Bigl\Vert_{p}.
\leqno(D'_{p})$$

It is easy to see that $(D_p)$ implies $(D'_{p'})$,
where $p'$ is the index conjugate to $p$. The other
implication is also easy in the commutative case.
However the non-commutative setting necessitates
more effort. We refer to [J] for more details. It
was also proved there that $\delta_{p'} =
\delta'_{p}$ for $1 \leq p < \infty.$

\medskip
In the rest of this paper, all the constants involved in the
preceding inequalities are assumed the best ones. Our aim is to
determine their optimal order of growth when $p\to 1$ or $\8$. We
will use the notation $a_{p} \approx b_{p}$ as $p \to p_{0}$ to
abbreviate the statement that there are two positive constants $c$
and $C$ such that
$$c \leq {a_{p} \over b_{p}} \leq C \quad \hbox{ for } p \hbox{ close to }
p_{0}.$$

Let us first recall the optimal order of these constants in the
commutative case.

\proclaim Optimal order in the commutative case. Only in this
statement that we use the same notations as before to denote the
best constants in the preceding inequalities in the {\it
commutative} case. Then their optimal orders are given as follows.
\itemitem{(i)} $\alpha_{p} \approx {(p-1)^{-1}}$ as $p \to 1$ \ ;\
$\alpha_{p} \approx \sqrt{p}$ as $p \to \infty.$
\itemitem{(ii)} $\beta_{p} \approx 1$ as $p \to 1$ ; $\beta_{p} \approx
p$ as $p \to \infty.$
\itemitem{(iii)} $\gamma_{p} \approx \sqrt{p}$ as $p \to \infty.$
\itemitem{(iv)} $\delta_{p} \approx {1 \over p-1}$ as $p \to 1.$
\itemitem{(v)} $\kappa_{p} \approx p$ as $p \to \infty.$

\medskip

We refer to [B1] for (i),(ii) and (v), to [B2] for (iv), and to
[St] for (iii). In the non-commutative setting, only very
recently, progress was made on these constants. It was proved in
[JX1] (using tools from [Mu]) that $\beta_{p}$ remains bounded as
$p \to 1$, as expected due to the fact that the second inequality
in $(BG_{p})$ remains true for $p = 1$ (cf. [PX1-2]). On the other
hand, Pisier [P] showed that $\beta_{p}={\rm O}(p)$ for even
integers $p$. As mentioned previously,  the spectacular progress
was achieved by Randrianantoanina [R]. He showed that $\kappa_{p}$
and $\beta_{p}$ have the same optimal order as in the commutative
case, and
$$\alpha_{p} \leq C p  \hbox{ as } p \to \infty,\quad \alpha_{p} \leq {C
\over (p-1)^{2}} \hbox{ as } p \to 1,\quad
\gamma_{p}\leq C p \hbox{ as } p \to \infty.$$
These  estimates directly follow from his  weak type
(1,1) estimate. Let us state this explicitly.

\proclaim Theorem 5 (Non-commutative martingale transforms: weak
type $(1,1)$). For all finite non-commutative $L^{1}$-martingales
$x=(x_n)$
$$\Bigl\Vert \sum_{n \geq 1} \varepsilon_{n} d x_{n}\Bigl\Vert_{1,
\infty} \leq C \Vert x \Vert_{1}, \ \forall\ \varepsilon_n = \pm
1, \leqno(M T_{1})$$
 where $\Vert\, \ldotp\, \Vert_{1, \infty}$ stands
for the non-commutative weak $L^{1}$-norm.

\medskip

Note that $(M T_{p})$ for $1 < p < \infty$ immediately follows
from $(M T_{1})$ by interpolation. This reduction of $(M T_{p})$
from $(M T_{1})$ yields $\kappa_{p} = {\rm O}(p)$ as $p \to
\infty,$ which is the optimal order. Randrianantoanina's proof for
$(M T_{1})$ heavily depends on a non-commutative version of the
classical Doob weak type $(1, 1)$ maximal inequality, obtained by
Cuculescu [Cu].

\medskip

Now we arrive at the position to state our main result. In order
to give a complete picture on the optimal orders of these
constants and to facilitate the comparaison with their commutative
counterparts as stated above, we incorporate in the following
statement the estimates  from [JX1] and [R].

\proclaim Theorem 6. We have the following estimates for the best
constants in $(BG_p)$, $(D_p)$, $(S_p)$ and $(MT_p)$.
\itemitem{(i)} $\alpha_{p} \approx p$ as $p \to \infty$; $\alpha_{p}
\leq {C (p-1)^{-2}}$ as $p \to 1.$
\itemitem{(ii)} $\beta_{p} \approx p$ as $p \to \infty$; $\beta_{p}
\approx 1$ as $p \to 1.$
\itemitem{(iii)} $\delta_{p} \approx (p-1)^{2}$ as $p \to 1.$
\itemitem{(iv)} $\gamma_{p} \approx p$ as $p \to \infty.$
\itemitem{(v)} $\kappa_{p} \approx p$ as $p \to \infty.$

\medskip

Thus the only problem left unsolved at the time of this writing is
on the optimal order of $\alpha_{p}$ as $p \to 1,$ which is
located between $(p-1)^{-1}$ and $(p-1)^{-2}$ (see the remark at
the end). As the reader can observe, compared with the commutative
case, these best constants can be divided into two groups,
according to their optimal orders. The constants in the first
group have the same optimal order as their commutative
counterparts; these constants are $\beta_{p}$ and $\kappa_{p}.$ On
the other hand, the optimal orders of those in the second one are
the squares of the respective optimal orders of their commutative
counterparts; the constants in this second group are $\alpha_{p}$
as $p \to \infty,$ $\delta_{p}$ and $\gamma_{p}.$

Our proof for the previous result mainly depends on the usual
triangular projection on the Schatten classes. Some elementary
facts on this projection and on the non-commutative martingales
with respect to  the natural filtration of matrices will be
presented in the next section. The third section will be devoted
to the proof of the result stated above. In the last section we
show in contrast to Davis' theorem in the commutative case that
the Hardy space ${\cal H}^1$ defined by  the square function does
not coincide with the Hardy space ${\cal H}^1_{\max}$ given by the
maximal function.

Let us end the present section by two further remarks. First, all
preceding inequalities still hold true in the more general
situation of Haagerup's non-commutative $L^{p}$-spaces (thus
including type III algebras), see [J] and [JX1]. We will show in a
work in preparation that all the best constants in this more
general setting remain the same as above. Second, in [JX2], we
proved the non-commutative ergodic maximal theorems. The same
approach based on Cuculescu's weak type (1,1) maximal inequality
yields an independent proof for the order $(p-1)^{-2}$ in Doob's
inequality.

\bigskip

\centerline{\bf 2. Canonical filtration of matrix algebras and
triangular projection}

\medskip

\n Recall that $S^{p}$ (resp. $S^{p}_{n})$ denotes the Schatten
$p$-class on $\ell_{2}$ (resp. $\ell^{n}_{2}).$ The elements in
these spaces as well as those in $B(\ell_{2})$ and
$B(\ell^{n}_{2})$ are represented as (finite or infinite)
matrices. For notational simplicity, we set $M_{n} =
B(\ell^{n}_{2})$ and $M_{\infty} = B(\ell_{2}).$ As usual, we
regard $M_{n}$ as a subalgebra of $M_{\infty}$ by viewing an $n
\times n$ matrix as an infinite one whose left upper corner of
size $n \times n$ is the given $n \times n$ matrix and all other
entries are zero. Note that $M_{n}$ is not a unital subalgebra of
$M_{\infty}.$ The unit of $M_{n}$ is the projection $e_{n} \in
M_{\infty}$ which projects a sequence in $\ell_{2}$ into its $n$
first coordinates.

Thus we have an increasing filtration $(M_{n})_{n \geq 1}$ of
subalgebras of $M_{\infty}$ whose union is w*-dense in
$M_{\infty}.$ This is the natural filtration of matrix algebras.
The corresponding conditional expectation from $M_{\infty}$ onto
$M_{n}$ is the mapping $E_{n}$ which leaves invariant the $n
\times n$ submatrix at the left upper corner of a matrix and
annihilates all other entries. Again, $E_{n}$ extends to a
contractive projection from $S^{p}$ onto $S^{p}_{n},$  still
denoted by $E_{n}.$ As in section 1, we define non commutative
martingales with respect to this canonical filtration $(M_{n})_{n
\geq 1}.$ We should emphasize that the present situation is
different from that in section 1 at two points. First, the
underlying von Neumann algebra $M_{\infty}$ is no longer finite
but semifinite; it is equipped with the usual normal semifinite
faithful trace ${\rm Tr}.$ Second, the conditional expectation
$E_{n}$ is no longer faithful; its support is the projection
$e_{n}$ above. In fact, we clearly have
$$E_{n} (a) = e_{n} a e_{n},\qquad a \in M_{\infty}.\leqno(1)$$
Neither of these two differences is essential for what follows.
Indeed, by approximation, we need only consider finite matrices,
i.e. those which have all but only finitely many zero entries;
then the trace ${\rm Tr},$ when restricted to each $M_{n},$ is
finite, and so can be normalized into a tracial state on $M_{n}.$
Thus in what follows, we can, and will,  ignore this difference.
As for as the second point, we need a little more effort. We have
to make all conditional expectations $E_{n}$ faithful, at least,
when restricted to finite matrices. This was already done in a
general setting in [JX1]. In the present setting the arguments in
[JX1] become much easier. Let us recall them briefly.

For a given finite matrix $a = (a_{ij}) \in M_{\infty}$ we define
$$\widetilde {E}_{n}  (a) = E_{n}(a) + \sum_{i > n} a_{ii}. \leqno(2)$$
Then $\widetilde {E}_{n}$ is faithful on the subalgebra of all
finite matrices. Thus when restricted to each $M_{n}$,
$\widetilde{E}_{1}, \cdots, \widetilde{E}_{n}$ form a finite
increasing filtration of faithful conditional expectations, and so
we are again in the situation described in section 1. Let $d_{1} =
E_{1},\ \widetilde{d}_{1} = \widetilde{E}_{1}$ and $d_{n} = E_{n}
- E_{n-1},\ \widetilde{d}_{n} = \widetilde{E}_{n} -
\widetilde{E}_{n-1}$ for $n \geq 2.$ Then for any $x \in S^{p},
(d_{n} x)_{n \geq 1}$ and $(\widetilde {d}_{n} x_{n})_{n \geq 1}$
are martingale difference sequences with respect to $(E_{n})_{n
\geq 1}$ and $(\widetilde{E}_{n})_{n \geq 1}$, respectively. We
have the following easily checked relations between these
martingale differences. Let $x = (x_{ij}) \in M_{\infty}$ be a
finite matrix. Then
$$\widetilde{d}_{1} x = d_{1} x + \sum_{i \geq 2} x_{ii}\quad \hbox{
and} \quad  \widetilde{d}_{n} x = d_{n} x -
x_{nn},\quad n \geq 2\ ;\leqno(3)$$

$${1 \over 3}\ \Bigl\Vert \Bigl( \sum_{n \geq 1}\ \vert d_{n} x
\vert^{2}\Bigl)^{1/2}\Bigl\Vert_{p} \leq \Bigl\Vert \Bigl( \sum_{n
\geq 1}\ \vert \widetilde{d}_{n} x\vert^{2}\Bigl)^{1/2}
\Bigl\Vert_{p} \leq 3 \ \Bigl\Vert \Bigl( \sum_{n \geq 1} \vert
d_{n} x \vert^{2} \Bigl)^{1/2} \Bigl\Vert_{p}. \leqno(4)$$ We
refer to [JX1] for the straightforward verifications and for more
information.

\medskip

Inequalities on non-commutative martingales with respect to
$(E_{n})_{n \geq 1}$ are closely related to the triangular projection
$T,$ which is defined by
$$(T(a))_{ij} = \cases{
a_{ij} &if $i \leq j$\cr \noalign{\medskip} 0 &otherwise\cr
},\quad a \in M_{\infty}.$$ It is classical that $T$ is bounded on
$S^{p}$ for $1 < p < \infty$ and the optimal order of the norm
$t_{p} = \Vert T \Vert_{S_{p} \to S_{p}}$ is ${\rm O}(p)$ as $p
\to \infty.$ Recall also that $T$ is selfadjoint on $S^{2},$ and
thus $t_{p} = t_{p'}.$ We will also need the norm of $T$ on
$S^{p}_{n}.$ Set
$$t_{p,n} = \Vert T \Vert_{S^{p}_{n} \to S^{p}_{n}}.$$
Thus $\displaystyle\sup_{n}\ t_{p,n} = t_{p}.$ It is also
classical that
$$t_{1, n} = t_{\infty, n} \approx\ \log\ (n+1), \quad \hbox{as}\ n\to\8.$$
We refer to [GK] and [KP] for the above classical facts on
$T.$

\bigskip

\centerline{\bf 3.  Proof of Theorem~6}
\medskip

\n In the following $\alpha_{p, n}, \beta_{p, n},$ $\ldots$ will
denote the best constants involved in the inequalities in
section~1 when restricted to all non-commutative martingales with
respect to the finite filtration $\widetilde{E}_{1}, \ldots,
\widetilde{E}_{n}$ on $S^{p}_{n}.$ By (1)-(3), we easily see that
these constants are equivalent uniformly on $p$ and $n$ to the
corresponding constants relative to all non-commutative
martingales with respect to the finite filtration $E_{1}, \ldots,
E_{n}$ on $S^{p}_{n}.$

\medskip

Now we proceed to the proof of Theorem~6. By the results from
[JX1] and [R] already quoted previously,  it remains to consider
$\alpha_{p}$ as $p \to \infty,\ \delta_{p}$ and $\gamma_{p}.$

\medskip

\proclaim Lemma 7. $\alpha_{\infty, n} \approx\ \log\ (n+1)$ as
$n\to\8$ and $\alpha_{p} \approx p$ as $p \to \infty.$

\n{\it Proof}. We consider the Hilbert matrix $\displaystyle h =
(h_{ij})_{1\le i,j\le n} \in M_{n}$ defined by
$$h_{ij} = \cases{
{(j-i)^{-1}} &if $i \not= j$\cr \noalign{\medskip} 0 &if $i =
j$\cr}.$$ It is well known that (cf. e.g. [KP])
 $$\Vert h \Vert_{\infty} \leq C \ \hbox{ and } \  \Vert T
 h\Vert_{\infty} \approx \log(n+1). \leqno(5)$$
For $1 \leq k \leq n$ we write $d_{k}h = a_{k} +
b_{k},$ where
 $$ a_{k} = (d_{k} h)(e_k-e_{k-1}) \quad \hbox{and}\quad
   b_{k} = (e_{k}-e_{k-1})(d_{k} h)$$
(recalling that $e_{k}$ is the natural projection on $\ell_{2}$
sending a sequence to its first $k$ coordinates). Note that $a_k$
and $b_k$ are the $k$-th column and $k$-th row of $d_k(h)$,
respectively. Moreover, we have  $a_{k} = d_{k} (T h).$ Thus
$a_{k}$ is the matrix whose $k$-th column is that of $T h$ and all
other columns are zero. Also note that
 $$\sum^{n}_{k=1} a_{k} a^{\ast}_{k} = (a_{1}, \ldots, a_{n})
 \pmatrix{a^{\ast}_{1}\cr \vdots\cr a^{\ast}_{n}\cr}.$$
It is
trivial  that the row matrix $(a_{1}, \ldots, a_{n})$ has the same
norm as $T h.$ Therefore, it follows that
$$\Bigl\Vert \Bigl( \sum_{k=1}^{n} a_{k} a^{\ast}_{k}\Bigl)^{1/2}
\Bigl\Vert_{\infty} = \Vert T h\Vert_{\infty}. \leqno(6)$$
On the other hand, $\displaystyle\sum^{n}_{k=1} b_{k} b^{\ast}_{k}$ is the
diagonal
matrix Diag $\Bigl( \displaystyle\sum^{k-1}_{k=1} j^{-2} \Bigl)_{1
\leq k \leq n}.$ Thus
$$\Bigl\Vert \Bigl( \sum^{n}_{k=1} b_{k}
b_{k}^{\ast}\Bigl)^{1/2}\Bigl\Vert_{\infty} = \Bigl(
\sum^{n-1}_{j=1} {1 \over j^{2}}\Bigl)^{1/2} \leq {\pi \over
\sqrt{6}}. \leqno(7)$$
Combining (5) - (7) we deduce
 $$\eqalign{
 \Bigl\Vert \Bigl( \sum^{n}_{k=1} \vert d_{k}(h)^{\ast}\vert^{2}
 \Bigl)^{1/2} \Bigl\Vert_{\infty} &\geq \Bigl\Vert \Bigl( \sum \vert
 a^{\ast}_{k}\vert^{2} \Bigl)^{1/2} \Bigl\Vert_{\infty} - \Bigl\Vert
 \Bigl( \sum \vert b^{\ast}_{k}\vert^{2}\Bigl)^{1/2}\Bigl\Vert_{\infty}\cr
 \noalign{\smallskip}
 &\geq C  \log (n+1) - {\pi \over \sqrt{6}}\cr
 \noalign{\smallskip}
 & \geq C'\log(n+1).\cr
 }$$
Then by (4) and (5), we get
$$\alpha_{\infty, n} \geq C \log(n+1).$$
The inverse inequality can be deduced from the known estimate
$\alpha_{p} \leq Cp$ as $p \to \infty$ obtained in [R]. Indeed,
for any $x \in S^{p}_{n}$
 $$
 \Bigl\Vert \Bigl( \sum^{n}_{k=1} \vert \widetilde{d}_{k} x
 \vert^{2} \Bigl)^{1/2}\Bigl\Vert_{\infty} \leq \Bigl\Vert \Bigl(
 \sum_{k} \vert \widetilde{d}_{k} x \vert^{2} \Bigl)^{1/2}
 \Bigl\Vert_{p} \leq \alpha_{p} \Vert x
 \Vert_{p} \leq \alpha_{p} n^{1/p} \Vert x \Vert_{\infty}\ ;
 $$
whence
$$\alpha_{\infty, n} \leq \alpha_{p}\ n^{1/p}.\leqno(8)$$
Choosing $p = \log (n+1)$ and using $\alpha_{p}
\leq C p$ as $p \to \infty,$ we get that
$\alpha_{\infty, n} \leq C \log (n+1).$ For the
second equivalence in the lemma it remains to prove
$\alpha_{p} \geq C p$ as $p \to \infty.$ But this
immediately follows from (8) by choosing this time
$n = [e^{p}]$. \qed

\smallskip

\n {\bf Remark.} This argument also shows that the lower
inequality in the Burkholder inequalities (see [JX1])  requires
the constant to be of order $p$. Indeed, we consider the Hilbert
matrix $h$ and $a_k$, $b_k$  as before. Note that
 $$ E_{k-1}(d_k(h)d_k(h)^*) = a_ka_k^* $$
and thus for this example
 $$ \| \sum_k
 E_{k-1}(d_k(h)d_k(h)^*)\|_{p/2}^{1/2} =
 \| T(h)\|_p \ge C p \ \|h\|_p \ .$$

\smallskip
\proclaim Lemma 8. (i) $\delta'_{p/2} \geq t^{2}_{p}$ for $2 \leq
p < \infty.$
\smallskip
\item{(ii)} $\gamma_{p} \geq t_{p}$ for $1 < p < \infty.$
\smallskip
\item{(iii)} $\gamma^{2}_{p} \leq \delta'_{p/2} \leq 2 \gamma^{2}_{p}$
for $2 \leq p < \infty.$

\n{\it Proof}. (i) Let $a \in M_{\infty}$ be a finite matrix. Let
$a_{k}$ be the matrix whose $k$-th column is that of $a$ and all
others are zero. Set $b_{k} =a_{k} a^{\ast}_{k}.$ Since $b_{k}
\geq 0,$ by (2)
$$E_{k} b_{k} \leq \widetilde{E}_{k} b_{k}, \quad k \geq 1.$$
Thus
$$\Vert \sum E_{k} b_{k}\Vert_{p/2} \leq \delta'_{p/2}\ \Vert \sum
b_{k}\Vert_{p/2}.$$
However,
$$\sum b_{k} = \sum a_{k} a^{\ast}_{k} = (a_{1}, a_{2}, \ldots)
\pmatrix{ a^{\ast}_{1}\cr a^{\ast}_{2}\cr \vdots\cr }.$$ As in the
proof of Lemma~7, the row matrix $(a_{1}, a_{2}, \ldots)$ has the
same norm as $a$ in $S^{p}.$ Thus
 $$\Vert \sum b_{k}\Vert_{p/2} = \Vert a \Vert^{2}_{p}\,.$$
On the other hand, by (1) and the same argument as above,
 $$
 \Vert \sum E_{k} b_{k}\Vert_{p/2} = \Vert \sum (e_{k} a_{k}) (e_{k}
 a_{k})^{\ast} \Vert_{p/2} =
 \Vert (e_{1} a_{1}, e_{2} a_{2}, \ldots )\Vert^{2}
 =
 \Vert T a \Vert^{2}_{p}.
 $$
Therefore, we deduce
$$\Vert T a\Vert^{2}_{p} \leq \delta'_{p/2} \Vert a \Vert^{2}_{p},$$
whence $t^{2}_{p} \leq \delta'_{p/2}$ for $2 \leq p \leq \infty.$

(ii) This can be proved in the same way as (i).

(iii) Let $\M$ and $({\cal E}_{k})_{k\ge 1}$ be as in section 1.
Let $(a_{k})$ be a finite sequence in $L^{p}({\cal M}).$ Since
$({\cal E}_{k} a_{k})^{\ast} ({\cal E}_{k} a_{k}) \leq {\cal
E}_{k} (a^{\ast}_{k} a_{k}),$ we have
$$\eqalign{
\Bigl\Vert \Bigl( \sum \vert {\cal E}_{k}
a_{k}\vert^{2} \Bigl)^{1/2} \Bigl\Vert_{p} &\leq
\Bigl\Vert \sum {\cal E}_{k} (a^{\ast}_{k} a_{k})
\Bigl\Vert^{1/2}_{p/2} \leq (\delta'_{p/2})^{1/2} \
\Vert \sum\ a^{\ast}_{k} a_{k} \Vert^{1/2}_{p/2}\cr
\noalign{\smallskip} &=(\delta'_{p/2})^{1/2}\ \Vert
(\sum a^{\ast}_{k} a_{k})^{1/2}\Vert_{p}.\cr }$$
Hence
$$\gamma_{p} \leq (\delta'_{p/2})^{1/2}.$$
The converse inequality is contained in [J]. \qed

\medskip

 Using Lemma~8 and the known estimate $\gamma_{p}
\leq C p$ as $p \to \infty$ from [R], we get the optimal orders of
$\delta_{p}$ and $\gamma_{p}$ in the theorem of section~1
(recalling that $\delta_{p'}=\delta'_{p}).$ Instead of using
$\gamma_{p} \leq C p$ from [R], we can use [JX2] to get
$\delta_{p} \leq C (p-1)^{-2}$ as $p \to 1$ ; and then by this and
Lemma~8 (iii) to recover $\gamma_{p} \leq C p$ as $p \to \infty.$
Therefore, we have completed the proof of Theorem~6.

\medskip\noindent {\bf Remark}.  The arguments in the proof of
Lemmas~7 and 8 show
$$\delta_{1,n} \approx (\log\ (n+1))^{2} \hbox{ and } \gamma_{1,
n} = \gamma_{\infty, n} \approx \log\ (n+1).$$

\bigskip

\centerline{\bf 4. Hardy spaces and maximal function}

\medskip\n Let us recall the classical Davis theorem (see [Da]) for
\underbar{commutative} martingales, namely
 $$ \| (\sum_k (dx_k)^2)^{1/2} \|_1 \sim_{c} \|\,
 \sup_k |{\cal E}_k(x)| \,  \|_1 \ . $$
If we denote by ${\cal H}^1_{\max}$ the space defined by the right
hand side, then this means that in the commutative case
 $$ {\cal H}^1 = {\cal H}^1_{\max}\,. $$
Using similar ideas as in the previous section, it turns out that
this equality does not hold in the non-commutative setting. For a
martingale sequence $(x_n)$ we  introduce the notation
 $$ \| (x_n)\|_{{\cal H}^p_{\max}} = \inf \| a\|_{2p}
 \ \sup_n \|y_n\|_\infty  \ \| b\|_{2p} \ <\ \infty
 \ , $$
where the infimum is taken over all $a,b\in L^{2p}({\cal M})$ and
bounded sequences $(y_n)\subset {\cal M}$ such that $x_n=ay_nb$
holds for all $n\in {\Bbb N}$ (see [J] for more information). Note
that for positive $(x_n)$ this is the same as
 $$ \| (x_n)\|_{{\cal H}^p_{\max}} = \inf\{ \|y\|_p
 \ : y\in L^p({\cal M}), \  \forall n\in {\Bbb
 N} \
  0\le x_n \le y  \}
 \ . $$
In the following lemma we refer to the faithful
filtration  $\widetilde {M}_{1},...,\widetilde
{M}_{n}$ defined by $\widetilde { M}_{k}=\widetilde
{E}_k(M_n)$ (see section 2)  for the algebra $M_n$
of $n\times n$ matrices.

\def\H{{\cal H}}

\medskip
\proclaim Lemma 9. $\|{\rm Id}:{\cal H}^1\to \H^1_{\max}\|\ge c
\log n$.

\n{\it Proof}.  Since the dual Doob inequality for $S_\8^n$ is
only valid with a constant $\|T\|^2=t_{1,n}^2$ (here $T$ denotes
the triangular projection on $S_1^n$), we deduce from Lemma 8 that
 $$ \|{\rm Id} : S_1^n\to \H^1_{\max}\|\ge t_{1,n}^2 \ .$$
It therefore suffices to show
 $$ \|{\rm Id} : S_1^n\to \H^1\|\le (1+2 t_{1,n})  \ .$$
Indeed, let $x\in S_1^n$, ${\rm Diag}(x)$ its diagonal and
$y=x-{\rm Diag}(x)$ the matrix with $0$ on the diagonal. Consider
$y_1=T(y)$ and $y_2=y-y_1$ obtained by applying the triangular
projection $T$. Then $d_k(y_1)$ are exactly the column matrices of
$y_1$ and
 $$ \| (\sum_k d_k(y_1)d_k(y_1)^*)^{1/2} \|_1 =\|
 y_1\| = \|T(x)\|_1 \le t_{1,n} \ \|x\|_1 \ .$$
Similarly,
 $$ \| (\sum_k d_k(y_2)^*d_k(y_2))^{1/2} \|_1 =\|
 y_2\| = \|T(x^*)\|_1 \le t_{1,n}\ \|x\|_1 \ .$$
Therefore,
 $$ \| x\|_{\H^1} \le \|{\rm Diag}(x)\| + 2t_{1,n} \ \|x\|_1 \le
 (1+2t_{1,n}) \ \|x\|_1 \ .$$
Thus
  $$ \|{\rm Id} : \H^1  \to \H^1_{\max}\| \ge
 {t_{1,n}^2\over 1+2t_{1,n}} \ .$$
$t_{1,n}\sim \log n$ implies that assertion.\qed

\proclaim Corollary 10. The spaces $\H^1$ and $\H^1_{\max}$ do in
general not coincide.

\medskip

\n {\bf Remark.} The idea in the proof above will be exploited
elsewhere in the study of $H^1$ and $BMO$ associated with a nest
algebra.

\medskip
\n {\bf Problem.} At the time of this writing the validity of the
inclusion
 $$ \H^1_{ \max}\subset \H^1 $$
is entirely open.

\medskip

\n{\bf Remark} (Added in March 2004). After this paper had been
submitted for publication, Randrianantoanina  finally settled up,
in January 2004,  the only case left unsolved in Theorem 6 on the
optimal order of $\alpha_{p}$ as $p \to 1$. He proved that
$\alpha_p\approx (p-1)^{-1}$ as $p\to 1$, the same optimal order
as in the commutative case.

\vfill\eject

 \centerline{\bf References}

\medskip

\item{[B1]} D. Burkholder, Sharp inequalities for martingales and
stochastic integrals, {\it Colloque Paul L\'evy sur les Processus
Stochastiques} (Palaiseau, 1987), {\it Ast\'erisque} No. 157-158
(1988), 75-94.

\item{[B2]} D. Burkholder, Exploration in martingale theory and its applications,
{\it Springer Lecture notes in Math.} 1464 (1991),
1-66.

\item{[CK]} E. Carlen, P. Kr\'ee, On martingale inequalities in non-commutative
stochastic analysis, {\it J. Funct. Anal.} 158 (1998), 475--508.

\item{[Cu]}I. Cuculescu,  Martingales on von Neumann algebras,
{\it J. Multivariate Anal.} 1 (1971), 17--27.

\item{[Da]} B. Davis, On the integrability of the martingale square
function, {\it Isreal J. Math.} 8 (1970), 187--190

\item{[D]} J. Dixmier, Formes lin\'eaires sur un anneau
d'op\'erateurs, {\it Bull. Soc. Math. France}, 81 (1953), 9--39.

\item{[GK]} I.C. Gohberg, M.G. Krein, Theory and applications of
Volterra operators in Hilbert space, {\it Transl. Math.
Monographs,} Vol. 24 {\it Amer. Math. Soc.} Providence, R.I. 1970.

\item{[J]} M. Junge,  Doob's inequality for non-commutative martingales,
to appear in {\it J. reine angew. Math.} 549 (2002), 149--190.

\item{[JX1]} M. Junge, Q. Xu, Non-commutative Burkholder/Rosenthal
inequalities, to appear in {\it Annals of
probability}.

\item{[JX2]} M. Junge, Q. Xu,  Maximal ergodic
theorems in non-commutative $L_p$ spaces, in preparation; see also
the note: Th\'eor\`emes ergodiques maximaux dans les espaces $L\sb
p$ non commutatifs, {\it C. R. Math. Acad. Sci. Paris} 334 (2002),
773--778.

\item{[KP]} S. Kwapie\'n, A. Pelczy\'nski, The main triangle
projection in matrix spaces and its applications, {\it Studia
Math.} 34 (1970), 43-67.

\item{[LPP]} F.\ Lust-Piquard, G.\ Pisier, Noncommutative Khintchine and Paley
inequalities, Arkiv f\"or Mat. 29 (1991), 241--260.

\item{[Mu]} M. Musat, Interpolation Between Non-Commutative BMO and Non-Commutative
$L_p$-Spaces, {\it J. Funct. Anal.} 202 (2003), 195--225.

\item{[NPTV]} F. Nazarov, G. Pisier, S. Treil, A. Volberg,
 Sharp estimates in vector Carleson imbedding theorem and
 for vector paraproducts, {\it J. Reine Angew. Math.}
 542 (2002), 147--171.

\item{[P]} G. Pisier, An inequality for $p$-orthogonal sums in non-commutative
$L_p$, {\it Illinois J. Math.}  44 (2000),  901--923.

\item{[PX1]} G. Pisier, Q. Xu, Non-commutative martingale inequalities,
 {\it Comm. Math. Phys.} 189 (1997), 667--698.

\item{[PX2]} G. Pisier, Q. Xu, In\'egalit\'es de martingales non commutatives,
{\it C. R. Acad. Sci. Paris} 323 (1996), 817--822.

\item{[PX3]} G. Pisier, Q. Xu, Non-commutative $L^p$-spaces,
Handbook of Geometry of Banach spaces, Vol. 2, pp. 1459-1517;
edited by W.B. Johnson and J. Lindenstrass, Elsevier Science,
2003.

\item{[R]} N. Randianantoanina, Non-commutative martingale
transforms,  {\it J. Funct. Anal.} 194 (2002), 181--212.

\item{[Ru]} M. Rudelson,  Random vectors in the isotropic position,
{\it J. Funct. Anal.} 164 (1999),  60--72.

\item{[S]} I.  Segal, A non-commutative extension of abstract
integration, {\it Ann. Math.} 57 (1953), 401--457.

 \item{[St]} E. M. Stein, Topics in harmonic analysis related to the
Littlewood-Paley theory, Princeton University Press, Princeton,
N.J., 1970.

\end